\documentclass[11pt]{amsproc}

\newtheorem{theorem}{Theorem}[section]
\newtheorem{lemma}[theorem]{Lemma}

\newtheorem{proposition}[theorem]{Proposition}
\newtheorem{corollary}[theorem]{Corollary}

\theoremstyle{definition}

\theoremstyle{remark}

\numberwithin{equation}{section}

\newcommand{\dach}[1]{\widehat{\vphantom{#1}}\;}

\def\G{\widehat{G}}
\newcommand{\dac}[1]{\widetilde{\vphantom{#1}}\;}
\def\mm{{\mathfrak m}}

\def \L{{\mathfrak L}}
\def \C{{\mathfrak C}}
\def \K{{\mathfrak K}}
\def\pext{{\rm Pext}}
\def\ext{{\rm Ext}}

\def\hom{{\rm Hom}}
\def\im{{\rm Im }\,}
\def\ker{{\rm Ker }\,}
\def\Z{{\bf Z}}
\def\Q{{\bf Q}}
\def\R{{\bf R}}
\def\T{{\bf T}}
\def \tildeG{\widetilde{G}}

\begin{document}

\title{Pure extensions of locally compact abelian groups}

\author{Peter Loth}
\address{Department of Mathematics\\
         Sacred Heart University\\
         5151 Park Avenue\\
         Fairfield, Connecticut 06825, USA}
\email{lothp@sacredheart.edu}


\begin{abstract}
In this paper, we study 
the group $\pext(C,A)$ for locally compact abelian (LCA) groups $A$ and $C$. 
Sufficient conditions are established for $\pext(C,A)$ to coincide with the
first Ulm subgroup of $\ext(C,A)$. Some structural information on pure injectives
in the category of LCA groups is obtained.
Letting $\C$ denote the class of LCA groups which can be written as
the topological direct sum of a compactly generated group and a discrete
group, we determine the groups $G$ in $\C$ which are pure injective in the category of
LCA groups. Finally we describe those groups $G$ in $\C$ such
that every pure extension of $G$ by a group in $\C$ splits and obtain a corresponding dual result.
\end{abstract}

\maketitle

\section{Introduction}
In this paper, all considered groups are Hausdorff topological abelian groups and will be 
written additively. 
Let ${\mathfrak L}$ denote the category of 
locally compact abelian (LCA) groups with continuous homomorphisms as morphisms. 
The Pontrjagin dual group of a group $G$ is denoted by $\G$ and the annihilator of
$S\subseteq G$ in $\G$ is denoted by $(\G,S)$. A morphism
is called {\em proper} if it is open onto its image, and a short exact sequence 
$$0\to A\stackrel{\phi}{\longrightarrow} B\stackrel{\psi}{\longrightarrow} C\to 0$$
in ${\mathfrak L}$ is said to be {\em proper exact} if $\phi$ and $\psi$
are proper morphisms. In this case, the sequence is called an 
{\em extension of $A$ by $C$ $($in $\L)$},
and $A$ may be identified with $\phi(A)$ and $C$ with $B/\phi(A)$. 
Following Fulp and Griffith \cite{FG1}, we let $\ext(C,A)$ denote the 
(discrete) group of extensions of $A$ by $C$. The elements represented
by pure extensions of $A$ by $C$ form a subgroup of $\ext(C,A)$ which is denoted by 
$\pext(C,A)$. This leads to a functor Pext from $\L\times\L$ into the category of discrete abelian
groups. The literature shows the importance of the notion of pure extensions (see for instance \cite{F}).
The concept of purity in the category of locally compact abelian groups 
has been studied by several authors (see e.g. \cite{A}, \cite{B}, \cite{Fulp1}, \cite{HH}, \cite{Kh},
\cite{L1}, \cite{L2} and \cite{V}). The notion of topological purity is due 
to Vilenkin \cite{V}: a subgroup $H$ of a group
$G$ is called {\em topologically pure} if $\overline{nH}=H\cap\overline{nG}$ for all
positive integers $n$. 
The annihilator of a closed pure subgroup of an LCA group is topologically pure (cf. \cite{L2}) but need not be pure 
in $\G$ (see e.g. \cite{A}). 
As is well known, $\pext(C,A)$ coincides with 
$$\ext(C,A)^1=\bigcap_{n=1}^\infty n\ext(C,A),$$
 the first Ulm subgroup of $\ext(C,A)$, provided that $A$ and $C$ are 
discrete abelian groups (see \cite{F}). 
In the category $\L$, a corresponding result need not hold: for groups $A$ and $C$ in $\L$,
$\ext(C,A)^1$ is a (possibly proper) subgroup
of $\pext(C,A)$, and it coincides with $\pext(C,A)$ 
if (a) $A$ and $C$ are compactly generated, or (b)
$A$ and $C$ have no small subgroups (see Theorem 2.4). 
If $G$ is pure injective in $\L$,
then $G$ has the form $R\oplus T\oplus G'$ where $R$ is a vector group, $T$ is a toral group
and $G'$ is a densely divisible topological torsion group having no nontrivial pure compact
open subgroups. However, the converse need not be true (cf. Theorem 2.7). Let
$\C$ denote the class of LCA groups which can be written as the topological direct sum
of a compactly
generated group and a discrete group. Then a group in $\C$ is pure 
injective in $\L$ if and only if it is injective in $\L$ (see Corollary 2.8). 
Let $G$ be a group in $\C$. Then every pure extension of $G$ by a group in $\C$ splits if and only if
$G$ has the form $R\oplus T\oplus A\oplus B$ where $R$ is a vector group, $T$ is a toral group,
$A$ is a topological direct product of finite cyclic groups and $B$ is a discrete bounded group.
Dually, every pure extension of a group in $\C$ by $G$ splits exactly if
$G$ has the form $R\oplus C\oplus D$ where $R$ is a vector group, $C$ is a compact torsion
group and $D$ is a discrete direct sum of cyclic groups (see Theorem 2.11).

The additive topological group of real numbers
is denoted by {\bf R}, {\bf Q} is the group of rationals, {\bf Z} is the group
of integers, {\bf T} is the quotient ${\bf R}/{\bf Z}$,
${\bf Z}(n)$ is the cyclic group of order $n$ and $\Z(p^\infty)$ denotes the quasicyclic group. 
By $G_d$ we mean the group $G$ with the discrete 
topology, $tG$ is the torsion part of $G$ and $bG$ is the subgroup of all 
compact elements of $G$. 
Throughout this paper the term ``isomorphic'' is used for ``topologically isomorphic'',
``direct summand'' for ``topological direct summand'' and ``direct product'' for
``topological direct product''. We follow the standard notation
in \cite{F} and \cite{HR}.

\section{Pure extensions of LCA groups}

We start with a result on pure extensions involving direct sums and direct products.
   
\begin{theorem}
Let $G$ be in $\L$ and suppose $\{H_i: i\in I\}$ is a collection of groups in $\L$. 
If $H_i$ is discrete for all but finitely many $i\in I$, then 
\begin{center}
$\pext(\bigoplus_{i\in I}H_i,G)\cong\prod_{i\in I}\pext(H_i,G).$
\end{center}
If $H_i$ is compact for all but finitely many $i\in I$, then 
\begin{center}
$\pext(G,\prod_{i\in I}H_i)\cong\prod_{i\in I}\pext(G,H_i)$.
\end{center}
In general, there is no monomorphism 
\begin{center}
$\pext(G,(\prod_{i\in I}H_i)_d)\rightarrow\prod_{i\in I}\pext(G,(H_i)_d).$ 
\end{center}

\end{theorem}

\begin{proof} 
To prove the first assertion, let $\pi_i: H_i\to\bigoplus H_i$ be the natural injection for 
each $i\in I$. Then the map  
$\phi:\ext(\bigoplus H_i,G)\to \prod\ext(H_i,G)$ defined by $E\mapsto (E\pi_i)$ is
an isomorphism (cf. \cite{FG1} Theorem 2.13), mapping the group $\pext(\bigoplus H_i,G)$ into 
$\prod\pext(H_i,G)$. If the groups $H_i$ and $G$ are stripped of their topology, the 
corresponding isomorphism maps the group $\pext(\bigoplus(H_i)_d,G_d)$ onto 
$\prod\pext((H_i)_d,G_d)$ (see \cite{F} Theorem 53.7 and p.$\;$231, Exercise 6). 
Since an extension equivalent to a pure extension is pure, $\phi$ maps $\pext(\bigoplus H_i,G)$ 
onto $\prod\pext(H_i,G)$, establishing the first statement. The proof of the second assertion is similar.
To prove the last statement, let $p$ be a prime and $H=\prod_{n=1}^\infty\Z(p^n)$, 
taken discrete. Assume $\ext(\widehat{\Q},H)=0$. By \cite{FG2} Corollary 2.10, the sequences
$$\ext(\widehat{\Q},H)\to\ext(\widehat{\Q},H/tH)\to 0$$
and
$$0=\hom((\Q /\Z)\dach{G},H/tH)\to\ext(\widehat{\Z},H/tH)\to\ext(\widehat{\Q},H/tH)$$
are exact, hence \cite{FG1} Proposition 2.17 yields $H/tH\cong\ext(\widehat{\Z},H/tH)=0$ which is
impossible. Since $\widehat{\Q}$ is torsion-free, it follows that $\pext(\widehat{\Q},H)=\ext(\widehat{\Q},H)\neq 0$.
On the other hand, we have
\begin{center}
$\prod_{n=1}^\infty\pext(\widehat{\Q},\Z(p^n))=\prod_{n=1}^\infty\ext(\widehat{\Q},\Z(p^n))\cong
\prod_{n=1}^\infty \ext(\Z(p^n),\Q)=0$
\end{center}
by \cite{FG1} Theorem 2.12 and \cite{F} Theorem 21.1. Note that this example shows that
Proposition 6 in \cite{Fulp1} is incorrect. 
\end{proof}

\begin{proposition}
Suppose $E_0: 0\to A\stackrel{\phi}{\longrightarrow}B\to C\to 0$ is a proper exact sequence
in $\L$. Let $\alpha : A\to A$ be a proper continuous homomorphism and $\alpha_{\ast}$
the induced endomorphism on $\ext(C,A)$ given by $\alpha_{\ast}(E)=\alpha E$. Then
$E_0\in\im\alpha_{\ast}$ if and only if $\im\phi/\im\phi\alpha$ is a direct summand of
$B/\im\phi\alpha$.
\end{proposition}

\begin{proof}
If $\alpha: A\to A$ is a proper morphism in $\L$, then
$$0\to\im\alpha\to A\to \im\phi/\im\phi\alpha\to 0$$
and
$$0\to\ker\alpha\to A\to\im\alpha\to 0$$
are proper exact sequences in $\L$ (cf. \cite{HR} Theorem 5.27). Now \cite{FG2} Corollary 2.10
and the proof of \cite{F} Theorem 53.1 show that $E_0\in\im\alpha_{\ast}$ if and only if the
induced proper exact sequence
$$0\to\im\phi/\im\phi\alpha\to B/\im\phi\alpha \to C\to 0$$
splits.
\end{proof}

If $A$ and $C$ are groups in $\L$, then $\ext(C,A)\cong\ext(\widehat{A},\widehat{C})$ (see \cite{FG1}
Theorem 2.12). We have, however:

\begin{lemma}
Let $A$ and $C$ be in $\L$. Then:
\begin{enumerate}
\item
In general, $\pext(C,A)\not\cong\pext(\widehat{A},\widehat{C})$.
\item
Let $\K$ denote a class of LCA groups satisfying the following property:
If $G\in\K$, then $\G\in\K$ and $nG$ is closed
in $G$ for all positive integers $n$. 
Then $\pext(C,A)\cong\pext(\widehat{A},\widehat{C})$ whenever $A$ and $C$ are in $\K$.
\end{enumerate}
\end{lemma}

\begin{proof}
(i) The finite torsion part of a group in $\L$ need not be a direct summand (see for instance
\cite{Kh}), so there is a finite group $F$ and a torsion-free group $C$ in $\L$ such that
$\pext(C,F)=\ext(C,F)\neq 0$. On the other hand, $\pext(\widehat{F},\widehat{C})\cong\pext(F,(\widehat{C})_d)=0$
by \cite{F} Theorem 30.2.

(ii) Let $A$ and $C$ be in $\K$ and consider the isomorphism 
$\ext(C,A)\stackrel{\sim}{\longrightarrow}\ext(\widehat{A},\widehat{C})$ given by $E: 0\to A\to B\to C\to 0 \mapsto 
\widehat{E}: 0\to\widehat{C}\to\widehat{B}\to\widehat{A}\to 0.$
The annihilator of a closed pure subgroup of $B$ is topologically pure in $\widehat{B}$ 
(cf. \cite{L2} Proposition 2.1) and
for all positive integers $n$, $nA$ and $n\widehat{C}$ are closed subgroups of $A$ and $\widehat{C}$,
respectively. Therefore, $E$ is pure if and only if $\widehat{E}$ is pure.
\end{proof}

Recall that a topological group is said to have {\em no small subgroups} if there is a
neighborhood of $0$ which contains no nontrivial subgroups. Moskowitz \cite{M}
proved that the LCA groups with no small subgroups have the form $\R^n\oplus \T^m\oplus D$ where
$n$ and $m$ are nonnegative integers and $D$ is a discrete group,
and that their Pontrjagin duals are precisely the compactly generated LCA groups.

\begin{theorem}
For groups $A$ and $C$ in $\L$, we have:
\begin{enumerate}
\item
$\pext(C,A)\supseteq\ext(C,A)^1$.
\item
$\pext(C,A)\neq\ext(C,A)^1$ in general.
\item
Suppose (a) $A$ and $C$ are compactly generated, or (b) $A$ and $C$ have no small subgroups.
Then $\pext(C,A)=\ext(C,A)^1$.
\end{enumerate}
\end{theorem}

\begin{proof}
(i) Let $\alpha: A\to A$ be the multiplication by a positive integer $n$ and let
$E: 0\to A\stackrel{\phi}{\longrightarrow} X\to C\to 0 \in n\ext(C,A)$. Since Ext
is an additive functor, there exists an extension $0\to A\to B\to C\to 0$ such that
$$
\begin{array}{ccrcccccc}
0 & \to & A & \to & B & \to & C & \to & 0\\
& & {\scriptstyle\alpha} \downarrow & & \downarrow & & \| & &\\
0 & \to & A & \stackrel{\phi}{\longrightarrow} & X & \to & C & \to & 0
\end{array}
$$
is a pushout diagram in $\L$. An easy calculation shows that $nX\cap\phi(A)=n\phi(A)$, hence 
$\ext(C,A)^1$ is a subset of $\pext(C,A)$.

(ii) Let $\pext(C,F)$ be as in the proof of Lemma 2.3. Then $\pext(C,F)\neq 0$ but 
$\ext(C,F)^1=0$.

(iii) Suppose first that $A$ and $C$ are compactly generated. If $\alpha:A\to A$ is the multiplication
by a positive integer $n$, then $\alpha(A)=nA$ is a group in $\L$. Since $A$ is $\sigma$-compact,
$\alpha$ is a proper morphism by \cite{HR} Theorem 5.29.
Let $E:0\to A\stackrel{\phi}{\longrightarrow} B\to C\to 0\in\ext(C,A)$.
By Proposition 2.2, $E\in\im\alpha_{\ast}=n\ext(C,A)$ if and only if
$\phi(A)/n\phi(A)$ is a direct summand of $B/n\phi(A)$. 
Now assume that $E$ is a pure extension. Then $\phi(A)/n\phi(A)$ is pure in the group $B/n\phi(A)$ which is
compactly generated (cf. \cite{M} Theorem 2.6). Since the compact group $\phi(A)/n\phi(A)$  
is topologically pure, it is a direct summand of $B/n\phi(A)$ (see \cite{L1} Theorem 3.1).
 Consequently, $E$ is an element of 
the first Ulm subgroup of $\ext(C,A)$ and by (i) the assertion follows. To prove the second part of (iii), assume that $A$ and $C$ have 
no small subgroups. By what we have just shown and Lemma 2.3, we have
$\pext(C,A)\cong\pext(\widehat{A},\widehat{C})=\ext(\widehat{A},\widehat{C})^1\cong\ext(C,A)^1.$
\end{proof}

By the structure theorem for locally compact abelian groups, any group $G$ in $\L$
can be written as
$G=V\oplus\tildeG$ where $V$ is a maximal vector subgroup of $G$ and $\tildeG$ contains a
compact open subgroup. The groups $V$ and $\tildeG$ are uniquely determined up to isomorphism
(see \cite{HR} Theorem 24.30 and \cite{AA} Corollary 1). 

\begin{lemma}
A group $G$ in $\L$ is torsion-free if and only if every compact open subgroup of $\tildeG$
is torsion-free.
\end{lemma}

\begin{proof}
Only sufficiency needs to be shown. Suppose every compact open subgroup of $\tildeG$ is torsion-free
and assume that $G$ is not torsion-free. Then $\tildeG$ contains a nonzero element $x$ of
finite order. If $K$ is any compact open subgroup of $\tildeG$, then $K+\langle x\rangle$ is
compact (see \cite{HR} Theorem 4.4) and open in $\tildeG$ but not torsion-free, a contradiction.
\end{proof}

Dually, we obtain the following fact which extends \cite{A} (4.33). Recall that a group is
said to be {\em densely divisible} if it possesses a dense divisible subgroup.

\begin{lemma}
A group $G$ in $\L$ is densely divisible if and only if $\tildeG/K$ is divisible for every
compact open subgroup $K$ of $\tildeG$.
\end{lemma}

\begin{proof}
Again, only sufficiency needs to be proved. Assume that $\tildeG/K$ is divisible for
every compact open subgroup $K$ of $\tildeG$ and let $C$ be a compact open subgroup of 
$(\G)\dac{G}$. Since $(\G)\dac{G}\cong(G/V)\dach{G}$ where $V$ is a maximal vector subgroup
of $G$, there exists a compact open subgroup $X/V$ of $G/V$ such that $C\cong ((G/V)\dach{G},X/V)
\cong((G/V)/(X/V))\dach{G}$ (see \cite{HR} Theorems 23.25, 24.10 and 24.11). By our
assumption, $(G/V)/(X/V)$ is divisible. But then $C$ is torsion-free (cf. \cite{HR} Theorem 24.23),
so by Lemma 2.5, $\G$ is torsion-free. Finally, \cite{R} Theorem 5.2 shows that $G$ is densely
divisible. 
\end{proof}

Let $G$ be in $\L$. Then $G$ is called {\em pure injective in $\L$} if for every pure extension
$0\to A\stackrel{\phi}{\longrightarrow}B\to C\to 0$ in $\L$ and continuous homomorphism
$f:A\to G$ there is a continuous homomorhism $\overline{f}:B\to G$ such that 
the diagram
$$
\begin{array}{ccrcccccc}
0 & \to & A & \stackrel{\phi}{\longrightarrow} & B & \to & C & \to & 0\\
& & {\scriptstyle f} \downarrow & \swarrow \,\scriptstyle\overline{f}& & & \\
& &  G & & & & 
\end{array}
$$
is commutative.
Following Robertson \cite{R}, we call $G$ a {\em topological torsion group} if 
$(n!)x\to 0$ for every $x\in G$. Note that a group $G$ in $\L$ is a topological torsion group if and only if
both $G$ and $\G$ are totally disconnected (cf. \cite{R} Theorem 3.15). Our next result
improves \cite{Fulp1} Proposition 9.

\begin{theorem}
Consider the following conditions for a group $G$ in $\L$:
\begin{enumerate}
\item
$G$ is pure injective in $\L$.
\item
$\pext(X,G)=0$ for all groups $X$ in $\L$.
\item
$G\cong \R^n\oplus \T^{\mm}\oplus G'$ where $n$ is a nonnegative integer, $\mm$ is a cardinal
and $G'$ is a densely divisible topological torsion group possessing no nontrivial pure
compact open subgroups.
\end{enumerate}
Then we have: $(i)\Leftrightarrow(ii)\Rightarrow(iii)$ and $(iii)\not\Rightarrow(ii)$.
\end{theorem}

\begin{proof}
If $G$ is pure injective in $\L$, then any pure extension $0\to G\to B\to X\to 0$ 
in $\L$ splits
because there is a commutative diagram 
$$
\begin{array}{ccccccccc}
0 & \to & G & \to & B & \to & X & \to & 0 \\
& & \| & \swarrow & & & & &  \\
& & \;G, & & & &
\end{array}
$$
hence (i) implies (ii).
Conversely, assume (ii). If $0\to A\to B\to X\to 0$ is a pure extension in $\L$ and $f:A\to G$
is a continuous homomorphism, then there is a pushout diagram
$$
\begin{array}{ccrcccccl}
0 & \to & A & \to & B & \to & X & \to & 0 \\
& & {\scriptstyle f} \downarrow & & \downarrow & & \| & & \\
0 & \to & G & \to & Y & \to & X & \to & 0.
\end{array}
$$
The bottom row is an extension in $\L$ (cf. \cite{FG1}) which is pure. By our assumption, it
splits and (i) follows.

To show (ii) $\Rightarrow$ (iii), let us assume first that $\pext(X,G)=0$ for all groups
$X\in\C$. Then the proof of \cite{L1} Theorem 4.3 shows that $G$ is isomorphic to
$\R^n\oplus \T^{\mm}\oplus G'$ where $n$ is a nonnegative integer, $\mm$ is a cardinal and
$G'$ is totally disconnected. Notice that $G'/bG'$ is discrete (cf. \cite{HR} (9.26)(a)) and
torsion-free. Since the sequence
$$0=\hom((\Q/\Z)\dach{G} ,G'/bG')\to\ext(\widehat{\Z},G'/bG')\to\ext(\widehat{\Q},G'/bG')=0$$
is exact, $G'/bG'$ is isomorphic to $\ext(\widehat{\Z},G'/bG')=0$ and therefore $G'=bG'$. It follows that
the dual group of $G'$ is totally disconnected (cf. \cite{HR} Theorem 24.17), thus $G'$ is
a topological torsion group. Suppose that $\pext(X,G)=0$ for all $X\in\L$ and let $K$ be a
compact open subgroup of $G'$. Then $G'/K$ is a divisible group (see \cite{Fulp2} Theorem 7 or the proof
of \cite{L1} Theorem 4.1), so by Lemma 2.6 $G'$ is densely divisible.
Now assume that $G'$ has a pure compact open subgroup $A$. Since $A$ is algebraically compact,
it is a direct summand of $G'$. But then $A$ is divisible, hence connected (see \cite{HR}
Theorem 24.25) and therefore $A=0$. Consequently, (ii) implies (iii).

Finally, (iii) $\not\Rightarrow$ (ii) because for instance, there is a nonsplitting extension
of $\Z(p^\infty)$ by a compact group (cf. \cite{A} Example 6.4).
\end{proof}
   
Those groups in $\C$ which are pure injective in $\L$ are completely determined:

\begin{corollary}
A group $G$ in $\C$ is pure injective in $\L$ if and only if $G\cong \R^n\oplus\T ^{\mm}$
where $n$ is a nonnegative integer and $\mm$ is a cardinal.
\end{corollary}

\begin{proof}
The assertion follows immediately from \cite{M} Theorem 3.2 and the above theorem.
\end{proof}

The following lemma will be needed.

\begin{lemma}
Every finite subset of a reduced torsion group $A$ can be embedded in a finite pure subgroup
of $A$.
\end{lemma}

\begin{proof}
By \cite{F} Theorem 8.4, it suffices to assume that $A$ is a reduced $p$-group. But then
the assertion follows from \cite{Ka} p.$\;$23, Lemma 9  and an easy induction.
\end{proof}

A pure extension $0\to A\to B\to C\to 0$ with discrete torsion group $A$ and compact group $C$ 
need not split, as
\cite{A} Example 6.4 illustrates. Our next result shows that no such example can occur if $A$ 
is reduced.

\begin{proposition}  
Suppose $A$ is a discrete reduced torsion group. Then $\pext(X,A)=0$ for all compactly generated groups
$X$ in $\L$.
\end{proposition}

\begin{proof}
Suppose $E: 0\to A\stackrel{\phi}{\longrightarrow} B\stackrel{\psi}{\longrightarrow} X\to 0$ 
represents an element of $\pext(X,A)$ where $A$ is a discrete reduced torsion group and $X$ is
a compactly generated group in $\L$. By \cite{FG2} Theorem 2.1, there is a compactly generated subgroup
$C$ of $B$ such that $\psi(C)=X$. If we set $A'=\phi(A)$, then $A'\cap C$ is discrete,
compactly generated and torsion, hence finite, so by Lemma 2.9 $A'$ has a finite pure subgroup $F$ 
containing $A'\cap C$. Now set $C'=C+F$. Then $F$ is a pure subgroup of $C'$ because it is pure
in $B$. But then $F$ is topologically pure in $C'$ since $C'$ is compactly generated. By
\cite{L1} Theorem 3.1, there is a closed subgroup $Y$ of $C'$ such that $C'=F\oplus Y$.
We have $B=A'+C=A'+C'=A'+Y$ and
$$A'\cap Y=C'\cap A'\cap Y= (F+C)\cap A'\cap Y=[F+(C\cap A')]\cap Y=F\cap Y=0,$$
thus $B$ is an algebraic direct sum of $A'$ and $Y$. Since $Y$ is compactly generated, it is
$\sigma$-compact, so by \cite{FG1} Corollary 3.2 we obtain $B=A'\oplus Y$.
Consequently, the extension $E$ splits.
\end{proof}

\begin{theorem}
Let $G$ be a group in $\C$. Then we have:
\begin{enumerate}
\item
 $\pext(X,G)=0$ for all $X\in \C$ if and only if
$G\cong\R^n\oplus\T^{\mm}\oplus A\oplus B$ where $n$ is a nonnegative integer, $\mm$
is a cardinal, $A$ is a direct product of finite cyclic groups and $B$ is a discrete
bounded group.
\item
$\pext(G,X)=0$ for all $X\in\C$ if and only if $G\cong\R^n\oplus C\oplus D$ where $n$
is a nonnegative integer, $C$ is
a compact torsion group and $D$ is a discrete direct sum of cyclic groups.
\end{enumerate}
\end{theorem}

\begin{proof}
Suppose $G\in\C$ and $\pext(X,G)=0$ for all $X\in\C$. By the proof of part (ii) $\Rightarrow$ (iii)
of Theorem 2.7, $G$ is isomorphic to $\R^n\oplus\T^{\mm}\oplus A\oplus B$ where $A$ is a compact totally
disconnected group and $B$ is a discrete torsion group. By Lemma 2.3, we have
$\pext(\widehat{A},X)\cong\pext(\widehat{X},A)=0$ for all discrete groups $X$, hence $\widehat{A}$ is
a direct sum of cyclic groups (see \cite{F} Theorem 30.2) and it follows that $A$ is a
direct product of finite cyclic groups. Again, we make use of \cite{A} Example 6.4 and conclude
that $B$ is reduced. But then $B$ is bounded since it is torsion and cotorsion. 
Conversely, suppose $G$ has the form $\R^n\oplus\T^{\mm}\oplus A\oplus B$ 
as in the theorem and let $X=\R^m\oplus Y\oplus Z$ where $Y$ is
a compact group and $Z$ is a discrete group. Then $\pext(X,A)\cong\pext(\widehat{A},\widehat{X})\cong\pext(\widehat{A},(\widehat{X})_d)=0$.
By Theorem 2.1, Proposition 2.10 and \cite{F} Theorem 27.5 we have 
$$\pext(X,B)\cong\pext(\R^m,B)\oplus\pext(Y,B)\oplus\pext(Z,B)=0$$
and conclude that
$$\pext(X,G)\cong\pext(X,\R^n\oplus\T^{\mm})\oplus\pext(X,A)\oplus\pext(X,B)=0.$$
Finally, the second assertion follows from Lemma 2.3 and duality.
\end{proof}

\end{document}